\documentclass{article}

\usepackage{arxiv}
\usepackage{tikz}
\usepackage[utf8]{inputenc} % allow utf-8 input
\usepackage[T1]{fontenc}    % use 8-bit T1 fonts
\usepackage{hyperref}       % hyperlinks
\usepackage{url}            % simple URL typesetting
\usepackage{booktabs}       % professional-quality tables
\usepackage{amsfonts}       % blackboard math symbols
\usepackage{nicefrac}       % compact symbols for 1/2, etc.
\usepackage{microtype}      % microtypography
\usepackage{lipsum}
\usepackage{graphicx}
\usepackage{siunitx}
\usepackage{float}
\usepackage{amsmath}
\graphicspath{ {./images/} }

\usepackage{tkz-euclide}
\usetikzlibrary{shapes,arrows}

\title{An Optimised Satellite Constellation for Forest Fire Detection through Edge Computing}

\author{
 Minduli Wijayatunga \\
  University of Sydney\\
  NSW, Australia \\
  \texttt{minduli1999@gmail.com} \\
   \And
 Prof Xiaofeng Wu \\
  University of Sydney\\
   NSW, Australia\\
  \texttt{xiaofeng.wu@sydney.edu.au} \\

}

\begin{document}
\maketitle
\begin{abstract}
The end of 2019 marked a bushfire crisis for Australia that affected more than $\SI{100000}{\kilo\meter\squared}$ of land and destroyed more than 2000 houses. In this paper, we propose a method of in-orbit fire detection to identify fires at their early stages, so that they can be eliminated before significant damage is done. An LEO satellite constellation is first developed through NSGA-II (Nondominated Sorting Genetic Algorithm II), optimising for coverage over Australia. Then, edge computing is adopted to run a bushfire detection algorithm on-orbit. A geostationary satellite is used for inter-satellite communication, as well as for maintaining a constant link to the ground so that bushfire detections can be reported back without a significant delay. Overall, this system is able to detect any fire that spans more than $\SI{5}{\meter}$ in length in under $\SI{1.39}{\second}$.
\end{abstract}

% keywords can be removed

\section{Introduction}

Australia has experienced bushfires for approximately 550 million years, starting with its separation from the Gondwana supercontinent \cite{Gond}. The average surface air temperature of Australia has anomolously increased over the past century due to climate change.  Hence, conditions that lead to forest fires such as high temperature, low humidity and low rainfall have been rising since the turn of the century \cite{cc}. As a result, Australia was in a bushfire crisis at the end of 2019. Within the state of New South Wales, the fires destroyed more than 2000 houses and forced many people to seek shelter elsewhere. By the end of 2019,$\SI{100000}{\kilo\meter\squared}$ of Australian land was scorched, causing 33 deaths \cite{b1}. \par 
As Australia has around 7.7 million km$^2$ of land area, it is difficult to detect small bushfires by ground-based observations alone. Failure to detect smaller bushfires can often result in increased spread and damage, as a bushfire can spread at a rate of  $\SI{2.5}{\meter\squared\per\second}$ \cite{firerate}. Currently, global coverage remote sensing satellites such as Landsat, Aqua, Terra and Sentinel are used to detect bushfires within Australia \cite{Terra}. These satellites are also used in the DEA hotspots service provided by Geoscience Australia \cite{dea}. However, as the revisit times of these satellites are very high, their observations are unable to generate real-time fire detections. The Himawari-8 satellite is another remote sensing satellite looking into bushfires. It can detect fires under 10 minutes as it is in a geostationary orbit \cite{him}. This satellite has channels of 2 km, 1km and 500 m resolutions, hence is only able to detect large scale fires\cite{him2}. There is no current method of fire detection that is both able to detect small fires and to make detections under 10 minutes. In this paper, we propose a constellation of remote sensing satellites for bushfire detection, equipped with a novel edge computing technology to achieve this capability. \par 
The following section describes the methodology we have used to create an edge computing satellite constellation. The overall results are then analysed in the following section, along with a discussion on how this project supplements the existing research on constellation development and edge computing. Lastly, a conclusion is given to summarise the results obtained.

\section{Methods}

%% Tikz flowchart
\tikzstyle{block} = [rectangle, draw, fill=white!20, 
    text width=10em, text centered, rounded corners, minimum height=4em]

\tikzstyle{arrow} = [thick,->,>=stealth]
\begin{figure}[h!]
\centering    
\begin{tikzpicture}[node distance = 4cm, auto]
    % Place nodes
    \node [block] (init) {Pick out image sensors for the constellation satellites.};
    \node [block, right of=init] (identify) {Calculate suitable orbits for the constellation satellites based on optimum coverage and transit time.};
    \node [block, right of=identify] (decide) {Develop a detection algorithm for bushfires.};
    \node [block, right of=decide] (B) {Implement edge computing to provide fast and efficient fire detections};
    % Draw edges
    \draw [arrow] (init) -- (identify);
    \draw [arrow] (identify) -- (decide);
    \draw [arrow] (decide) -- (B);
\end{tikzpicture}
\caption{Approach taken for fast and efficient fire detection}
    \label{kelp}
\end{figure}
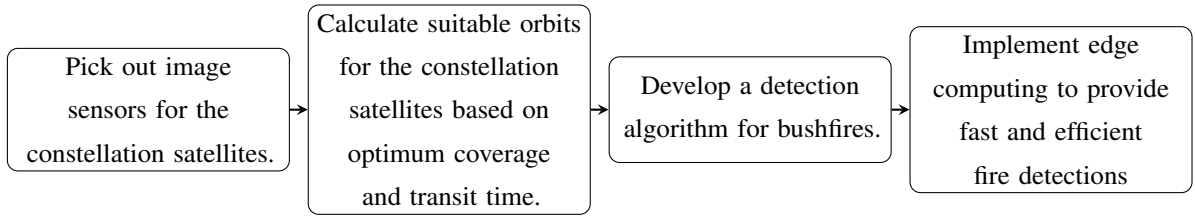
Figure \ref{kelp} shows the steps taken to achieve fast and efficient fire detections in this project.

\subsection{Choosing Imaging Sensors for the Satellite Network}
 Being able to detect small scale fires require high-resolution imaging \cite{modis1}. However, high resolutions lead to lower fields of view, which increase the number of satellites required to observe the area of interest \cite{twelve}. Table \ref{t7} describes several imaging sensors that are currently in use by LEO satellites. \par 
 
 \begin{table}[H] 
\caption{Imaging payloads that are currently used on-board satellites}
\label{t7}
%%% \tablesize{} %% You can specify the fontsize here, e.g., \tablesize{\footnotesize}. If commented out \small will be used.
\begin{tabular}{lc}
\toprule
\textbf{Imaging payload}	& \textbf{resolution}\\
\midrule
RapidEye Constellation sensors \cite{RE}                  &$\SI{5}{\meter}$\\
Flock Constellation sensors \cite{five}  &$\SI{3}{\meter}$ - $\SI{5}{\meter}$\\
Landsat-8 sensor \cite{five}  &$\SI{30}{\meter}$\\
Disaster monitoring Constellation (DMC) sensors \cite{DMC}  &$\SI{22}{\meter}$\\
Spire Lemur Constellation (DMC) sensors \cite{SLC}  &$\SI{1}{\kilo\meter}$\\
\bottomrule
\end{tabular}
\end{table}

For our project, we propose using sensors that have properties similar to Planet's RapidEye sensors. These sensors have multi-spectral detection capabilities in the visible and infrared range, along with a field of view of 13.08 degrees and a pixel size of 5m \cite{re1}. These sensors can give us a sufficiently high image resolution for fire detection without compromising the field of view drastically. \par

\subsection{Calculating Suitable Orbits for the Constellation Satellites}
In this section, we optimise the coverage and revisit times of the constellation to find the best orbits for its satellites. We use NSGA-II (Non-dominated Sorting Genetic Algorithm II) for this process as it is faster to run compared to other optimisation methods and seldom gets stuck in local optima \cite{theiss}. 

\subsubsection{Structure of NSGA II}\label{GA}
NSGA-II consists of several distinct components that are continuously reused. These are the chromosome encoding, the fitness function, selection, recombination and the scheme of evolution \cite{nine}. 
\begin{itemize}
    \item \textbf{Chromosome encoding}-
A chromosome in NSGA-II represents a single possible solution to a particular problem. NSGA-II uses a string of bits to represent a chromosome. As some problems being solved by genetic algorithms can contain large solutions, a GA can operate on bit strings of lengths greater than 100 \cite{nine}.
\item \textbf{The fitness function}-
By definition, the fitness function measures the suitability of a particular chromosome as the solution to the given problem \cite{nine}. Suppose the problem is to minimise or maximise a certain quantity (for example, minimising the revisit time of our satellite constellation). In that case, the fitness function will provide us with a measure of that quantity (sum of revisit times for all satellites in the constellation).
Problem-specific information is entirely contained in the fitness evaluation. Thus, the same genetic algorithm can be implemented for two entirely different problems just by changing the fitness function. 
\item \textbf{Selection} - Selection is the process in which the genetic algorithm uses the fitness function as a guide to offer selective pressure to the chromosomes under evolution. This process causes the chromosomes to be selected based on fitness for recombination. The chromosomes of higher fitness are selected more often than those of low fitness. Selection happens without replacement so that chromosomes of high fitness can be selected several times by the genetic algorithm \cite{nine}. \par 
Generally, the probability of a chromosome being selected is promotional to its relative fitness. However, other, more complex selection schemes such as Tournament Selection and Truncation Selection are also commonly used in genetic algorithms \cite{nine}.
\item \textbf{Recombination}-
Recombination is the process in which the selected chromosomes are reconfigured to form new members of the succeeding population. The recombination process mimics the process of reproduction described in Darwin's theory of evolution. \par 
Recombination happens through both crossover and mutation. Both are non-deterministic, as each occurs with a given probability that can be set by the user. The process of crossover represents the mixing of the genetic material of two selected chromosomes to produce a 'child' chromosome \cite{nine}.  An example of a crossover in progress is shown in Figure \ref{crossover}.

\begin{figure}[h!]
    \centering
    \includegraphics[width = 0.3\textwidth]{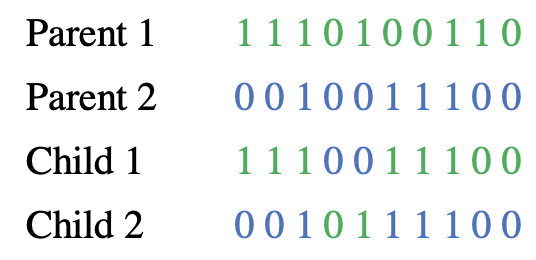}
    \caption{Crossover between two parent chromosomes to create a child chromosome}
    \label{crossover}
\end{figure}

In Figure \ref{crossover}, child 1 takes the first 4 bits from parent one and the last 4 bits from parent two, while child 2 does the reverse. \par
After crossover, the child chromosomes are passed on to the mutation operators. They act on the child chromosomes and flip one or more individual bits.  
\item \textbf{Evolution}-
Once the process of recombination has been completed, the child chromosomes are added to the succeeding population. Then the selection and recombination processes are repeated on all individual chromosomes of the parent population till a complete succeeding population is formed. Then the succeeding population is denoted as the current population and undergoes the same process to generate another population (next generation) \cite{nine}. \par 
The genetic algorithm can go through several generations to reach the required criteria. Other stopping criteria such as a maximum number of generations are often introduced to ensure that the GA does not run forever \cite{nine}. 
\end{itemize}

\subsubsection{Applying NSGA-II for Calculating Constellation Orbits}
Only the chromosome and fitness function vary from problem to problem when implementing NSGA-II. Thus, we have to wisely choose both before the start of the optimisation process. \par 
\textbf{Chromosome encoding}\par 
After analysing the outcomes of representing the constellation as a whole and as individual satellites in the optimisation process, it can be found that representing the constellation as a whole can generate orbits with more coverage using less total number of satellites \cite{theiss}. Thus, we use this technique to generate a chromosome. A chromosome (x) is then encoded to contain six genes in the order given in Table \ref{pf2}.\par 

\begin{table}[H] 
\caption{Constellation information allocation in a chromosome}
\label{pf2}
%%% \tablesize{} %% You can specify the fontsize here, e.g., \tablesize{\footnotesize}. If commented out \small will be used.
\begin{tabular}{cc}
\toprule
\textbf{Parameter}	& \textbf{Position in the chromosome}\\
\midrule
Semi major axis ($a$)                   & $x(1)$ \\
Eccentricity ($e$)                     & $ x(2)$ \\
Inclination ($i$)                     & $x(3) $ \\
Number of orbital planes (P) & $ x(4)$ \\
Relative spacing between satellites in adjacent planes (F) & $x(5) $ \\
Number of satellites per plane (n)                 &  $x(6)$\\
\bottomrule
\end{tabular}
\end{table}

This method of constellation representation is often used to represent Walker constellations \cite{walker2}. Ascending nodes of each plane of satellites is placed $360/P$. Thus, the right ascension of the ascending nodes ($\Omega$) of each plane is $360/P$ degrees apart. In each plane, the satellites are separated $360/n$ degrees. Thus, the argument of periapsis ($\omega$) of satellites in a single plane is $360/n$  degrees apart. Lastly, the true anomaly of equivalent satellites in adjacent planes is $360F/N$ degrees apart, where N is the total number of satellites, obtained by $N = nP$. Thus, the mean anomalies are also  $360F/N$ degrees apart \cite{Walker}. These relations allow us to convert these six parameters to Keplerian parameters. \par 
The upper and lower limits of the genes in our chromosome also need to be set so that our solution space is in the expected region of space. We are searching for a low Earth orbit constellation; thus, we restrict the orbital altitudes to 200 km - 1000 km. Thus
$ (R_e +200)\SI{}{\kilo\meter}  \leq a  \leq (R_e +1000)\SI{}{\kilo\meter} $, where $R_e$ is the radius of Earth and is taken to be $\SI{6371}{\kilo\meter} $ \cite{Earth}. As $a$ is relatively small, we cannot make our orbits highly eccentric without having their perigees be inside (or on the surface of) the planet. Thus we restrict our range of eccentricities to  $0 \leq e \leq 0.05$. We also limit P to 1-100, n to 1-50 and F to 1-8 such that a constellation with at most 5000 satellites can be devised. It can be noted that as we can optimise P and n, we can optimise the total number of satellites in the constellation.  \par 

\textbf{Fitness Function Determination}\par 
Now that we have our chromosomes encoded and the limits of each gene set, we go on to the derivation of the fitness function. Due to the nature of optimisation algorithms, the output of the fitness function must be a quantity to be minimised \cite{GA}. In our case, we need a fitness function output that is proportional to the total number of satellites and inversely proportional to the coverage over Australia. \par 

A variable that represents the coverage is the number of satellites above Australia at a given time. We shall call this variable C from now on. The more satellites that are over the country, the greater proportion of the land is observed.  \par 

\begin{figure}[h!]
    \centering
  \includegraphics[width=0.9\textwidth]{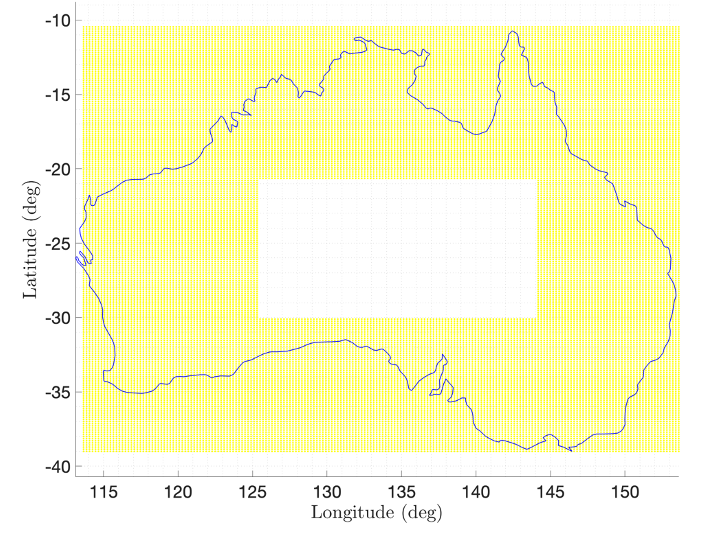}
  \caption{Area expected to be covered by the constellation within the Australian mainland. Each yellow point denotes a gridpoint.}\label{Australia}
\end{figure}

Another variable that also represents coverage is obtained by the process of constructing a  $200 \times 200$ grid through our area of interest (shown in Figure  \ref{Australia}). The latitude and longitude of each grid point are determined, and the number of grid points that the constellation covers within a day is counted. We call this parameter P (for points covered). As the middle of Australia does not have the conditions to grow forests, they seldom experience bushfires \cite{b1}. Thus, we exclude the grid points in the middle when calculating P. 

\begin{figure}[h!]
  \centering
  \includegraphics[width=0.8\textwidth]{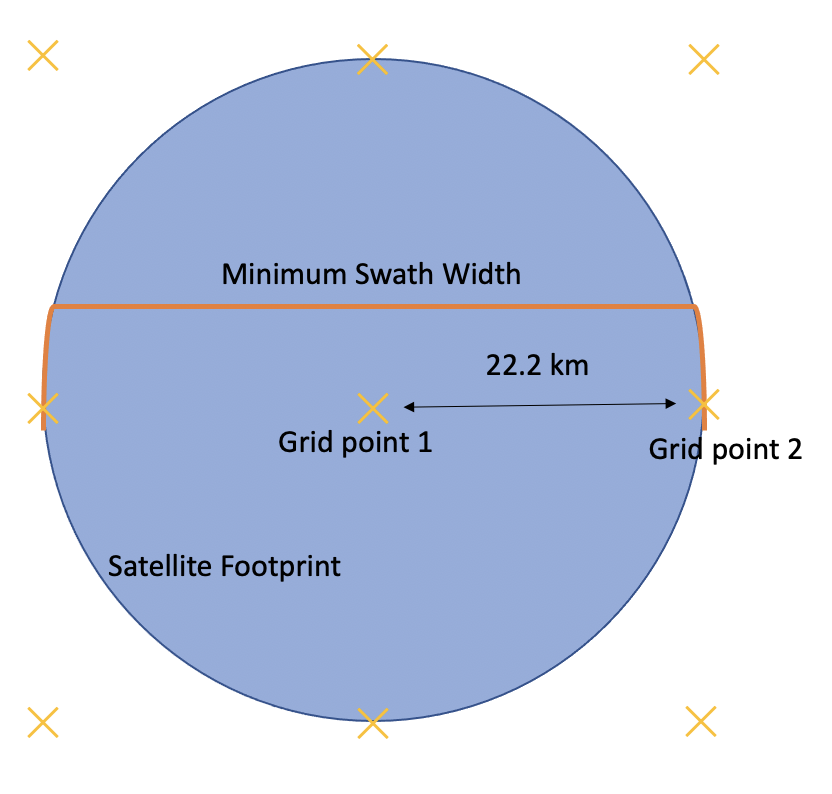}
  \caption{The distance between two grids being 22.2 km requires each satellite to have a swath width of at least 44.4 km ($2 \times 22.2$ km) to ensure full visibility when all grid points (shown as yellow crosses) are covered.}\label{explain}
\end{figure}

Note that the $200 \times 200$ grid constructed by the fitness function allows the distance between any two grid points to be $ 22.2$ km \cite{distance}. Thus, when all grid points are covered, there would be 100\% coverage as long as the swath widths of the constellation satellites are greater than $44.4$ km, as shown in Figure \ref{explain}. Thus, we introduce this condition to the fitness function. \par

\bigskip
A variable that is proportional to the revisit time is the average time between two satellite transits across the area of interest. We shall call this variable R. This does not represent the revisit time of an individual satellite, but the revisit time of the constellation.  \par
We can expect that the coverage and revisit times would increase with the number of satellites \cite{theiss}. However, as the launch and construction of satellites is an expensive process, we require the minimum number of satellites that still gives full coverage of the area of interest. Therefore, we must minimise the total number of satellites in the constellation (N). \par 
Thus, the fitness function goes through the following steps to determine C, P and R and to combine them with N to form an output.

\begin{enumerate}
\item Create a matrix of instances($I$) to store the position data created. 
\item Create a variable P to store the number of grid points crossed as each satellite passes above the area of interest.  
    \item For satellite number $i = 1$ to $i = N$, do the following. 
\begin{enumerate}
    \item \textbf{Obtain the Earth-centred inertial (ECI) coordinates of the satellite for 1 Earth Day.}\par 
    
First, we calculate the mean anomalies (MA) within the satellite orbit using Equation \ref{MA}. 
\begin{equation}
     \text{MA} = \text{MA} + \sqrt{\frac{\mu}{a^3}}(t- t_0)
     \label{MA}
\end{equation}
Here, $\mu$ is the gravitational constant and t is the time. We set $t_0$ (time at perigee) to be equal to 0, for ease of optimisation. 
Then we use the Newton-Raphson method to obtain the Eccentric anomalies(E) from the Mean anomalies calculated, using Equation \ref{E}. 
\begin{equation}
     MA = E - e \sin{E}
     \label{E}
\end{equation}
The true anomaly (f) is then calculated using Equation \ref{f}.
\begin{equation}
     f = 2 \tan^{-1}\left(\sqrt{\frac{1 +e}{1 -e}}\tan\frac{E}{2}\right)
     \label{f}
\end{equation}
The radial distance to the satellite from the centre of the Earth at each point of its orbit is then calculated. 
\begin{equation}
    r = \frac{a(1-e^2)}{1 + e \cos{f}}
\end{equation}
Then the 2-D orbit coordinates are determined.
\begin{equation}
    \begin{bmatrix}X \\ Y \\ Z
    \end{bmatrix}
    =
    \begin{bmatrix} R\cos{\theta}\\ R\sin{\theta}\\  0
    \end{bmatrix}
\end{equation}
    
    Finally, the ECEF coordinates are determined by multiplying the 2-D coordinates by the following conversion matrix. 
    
    \begin{equation}
         C^{ECI}_{Orbit} = 
\begin{pmatrix}
 \cos{\omega}\cos{\Omega}-\sin{\Omega}\sin{\omega}\cos{i} &  -\cos{\Omega}\sin{\omega}-\sin{\Omega}\cos{\omega}\cos{i}&  \sin{\Omega}\sin{i}
 \\ 
  \sin{\Omega}\cos{\omega} + \cos{\Omega}\sin{\omega} \cos{i} & 
  -\sin{\omega}\sin{\Omega} + \cos{\Omega}\cos{\omega}\cos{i}& 
  -\cos{\Omega}\sin{i}
  \\
  \sin{\omega}\sin{i} & 
  \cos{\omega}\sin{i} & 
  \cos{i} \end{pmatrix}
    \end{equation}

    \item \textbf{Obtain the latitudes and longitudes traversed by the satellites through coordinate conversions.}\par 
    First, we have to convert the ECI coordinates to Earth-centred Earth Fixed(ECEF) coordinates by multiplying them by the following conversion matrix. 

Given that:
$$\omega = \SI{7.292e-5}{} = \text{Frequency corresponding to a sidereal day}$$
$$ t = \text{Time since the vernal equinox}$$
$$
 C^{ECEF}_{ECI} = 
 \begin{bmatrix}
 \cos{\omega t} & \sin{\omega t} & 0 \\
 -\sin{\omega t } & \cos{\omega t} & 0 \\ 
 0 & 0 & 1 
  \end{bmatrix}
$$
Then we convert the ECEF coordinates to longitude as follows. 
\begin{equation}
    \text{Longitude} = \tan^{-1}\left(\frac{y_{ECEF}}{z_{ECEF}}\right)
\end{equation}
The latitude is then determined iteratively using Equation \ref{lat}, where $R_E$ is the radius of Earth. 
\begin{equation}
    \text{Latitude} = \tan^{-1}\left(\frac{z + e^2 R_{E} (\frac{z}{R_{E}( 1 - e^2)} + h )}{\sqrt{x^2 + y^2}} \right)
    \label{lat}
\end{equation}
\item If at time $t$, the latitude and longitude lies within the yellow section in Figure \ref{Australia}, denote: $$ \text{instance(t)} = 1 $$
 \item 
    Count the number of grid points crossed by the satellite once it is in the yellow region in Figure \ref{Australia}. Add this to P. 
  
    \item If at time $t$, it does not lie within the yellow region, denote: $$ \text{instance(t)} = 0 $$
    \item Now, for satellite(i), a vector of instances is generated. This vector has 1 for instances when the satellite can see the chosen area of interest, and 0 when that satellite cannot. Let's call this vector: $$\text{instances}_i$$
      \end{enumerate}
          \item The final $I$ matrix should take the following form. 
     \begin{equation}
          I = \begin{bmatrix}  \text{instances}_1 & \text{instances}_2 & \hdots & \text{instances}_i & \hdots & \text{instances}_{N} \end{bmatrix}^T
      \end{equation}
      \item To determine C, sum the instances matrix along its columns. This tells us how many satellites are visible at each time step. Taking the average of this gives us the C parameter. 
      \begin{equation}
          C = avg( \sum^{j = j_{max}}_{1} I_{ij})
      \end{equation}
          \item To determine R, calculate the space between the consecutive 1s in the I matrix, and find the average. This gives the average time gap between transits.  
          \item We have to maximise C and P and minimise R and N. Thus, the fitness function used is the following. As all parameters are of equal importance, we need to ensure that they are all on the same scale. Thus, we introduce a factor of 0.01 in front of R.
          \begin{equation}
              \text{Fitness function} = \frac{1}{C} + 0.01R + \frac{1}{P} \ + N
          \end{equation}
\end{enumerate}

\subsection{Developing a Detection Algorithm for Bushfires}
The second part of this project focuses on the choosing an on-orbit image processing method to detect fires.  Over the years, several different algorithms have been constructed for satellite-based image processing. Some of the methods are as follows. 
\begin{enumerate}
    \item \textbf{Collection 6 MODIS active fire detection algorithm \cite{C6MODIS}}\label{c6}\par 
    This algorithm is based on the data collected by the Moderate Resolution Imaging Spectroradiometer (MODIS), onboard the NASA satellites, Terra and Aqua. Data is obtained from 7 different channels with variations in central wavelength, and are then used to determine if a certain pixel depicts fire or not. The algorithm involves water and cloud masking for efficiency, accompanied by sun glint rejection for improved accuracy \cite{C6MODIS}. \par 
This algorithm is useful to us as it is implemented on 2 LEO satellites, both within the range of altitudes that we are concerned with. Furthermore, the error associated with a false detection using this method is derived to be 1.2\% \cite{C6MODIS}, confirming that the Collection 6 algorithm is an adequate fire detector. 
    \item \textbf{The enhanced contextual fire detection method \cite{EMODIS}}\par 
    This algorithm has a similar basis to the previous and is also based on MODIS data. However, it offers better sensitivity to smaller fires as well as a low chance of false alarms \cite{EMODIS}. The  250-m resolution red and near-infrared channels are used in this algorithm, for the rejection of false detections. 
    \item \textbf{Fire detection using Infrared image sensors  \cite{MSG}}\par 
     In his paper, Giovanni Laneve focuses on using an Infrared sensor onboard a geostationary satellite for fire detections in the mediterranean region \cite{MSG}.
     He notes that a satellite-based fire detection is useful if it can detect fires at least 1500 m$^2$ small with a 30 min temporal resolution \cite{MSG}. \par 
The algorithm proposed uses a change detection technique that compares two or more images obtained at 15-minute intervals. If we were to run this algorithm onboard a satellite, previously taken images would need to be sent from one satellite to another. However, as we discussed in the Introduction section, this process is too time-consuming to be feasible for useful bushfire detection. Thus, this algorithm cannot be used onboard a satellite.
    
    \item \textbf{Machine learning-based fire detection \cite{ml}} \par 
    The article titled "Detection and Monitoring of Forest Fires Using Himawari-8 Geostationary Satellite Data in South Korea" by Eunna Jang develops a fire detection algorithm with three steps, namely, thresholding, machine learning and post-processing. The algorithm is implemented on Himawari-8 geostationary satellite data. Thresholding is used to filter candidate pixels of a forest fire, while random forests is used to remove any non-fire pixels filtered in when thresholding. Under this algorithm, the probability of false detection drops down to 0.07\% \cite{ml}. However, this algorithm is reliant on past data.
\end{enumerate}

\subsubsection{Algorithm Chosen for Onboard Fire Detection}\label{cover}
As seen in the previous section, The Collection 6 bushfire detection algorithm is proven to provide the most accurate fire detections from images of Australia taken by the Terra satellite during the bushfire season \cite{theiss}. The detailed algorithm is given in \cite{modis1}. However, this algorithm is intended to be used for ground-based processing, as parts of it rely on image comparisons.  In this project, we adapt the Collection 6 detection algorithm for onboard fire detection by modifying the parts of it that depend on image comparisons \cite{theiss}. The Collection 6 fire detection method goes through each pixel in an image and recognises whether if it depicts a fire or not. \par 
The inputs required for the algorithm are as follows. 
\begin{itemize}
    \item 
\textbf{Detection channels of 5 m resolution of varying central wavelengths in the infrared range} \par 
The Collection 6 fire detection algorithm requires observations obtained using channels at  wavelengths  0.65$\SI{}{\micro\meter}$ , 0.86$\SI{}{\micro\meter}$ , 4.0$\SI{}{\micro\meter}$ , 11.0$\SI{}{\micro\meter}$  and 12.0$\SI{}{\micro\meter}$. The first two channels are only required to make detentions in the daytime, while the latter three are required for day and night times \cite{modis1}. These are all in the long-wave thermal infrared range \cite{TIR}. 
\item \textbf{Solar zenith angle at the time of observation}\par 
The solar zenith angle is the angle between the Sun's position and the local zenith  (i.e. the point directly above the ground) \cite{sza}. It is required to determine whether it is day or night time at the location being imaged. The standard is to assume daytime if the solar zenith angle is below 85 degrees \cite{modis1}.

\item \textbf{Water mask (preferred)}\par 
A water mask is preferred for the algorithm, as removing water pixels speeds up the fire detection process. This mask can be generated using the SRTM Water Body Dataset \cite{SRTM}.
\end{itemize}

\bigskip

The steps of the algorithm are as follows. 
\begin{enumerate}
    \item 
    The Plank's law is used to convert the radiance of each pixel to a brightness temperature as the brightness temperatures of radiances obtained in 4.0$\SI{}{\micro\meter}$, 11.0$\SI{}{\micro\meter}$  and 12.0$\SI{}{\micro\meter}$ channels are required for the following calculations \cite{b1}. 
\begin{equation}
    T = \frac{c_2}{L \ln(\frac{c_1}{b L^5} +1)} \ \text{where} \  c_1 = 2 h c^2, \ \text{and} \  c_1 = \frac{hc}{k}
\end{equation}
Here, $T,L,b,h,c$ and $k$ represent the brightness temperature in Kelvin, the central wavelength in metres, reflectances in Watts per square meter per steradian per micron , Plank's constant, speed of light, and the Boltzmann's constant, respectively.
\item 
These brightness temperatures are then denoted as $T_4, T_{11}$ and $T_{12}$. The reflectance values obtained from the 0.65$\SI{}{\micro\meter}$ and 0.86$\SI{}{\micro\meter}$  channels are denoted as $\rho_{0.65}$ and $\rho_{0.86}$, respectively.
\item 
Any pixels that are not on land are removed using the water mask.
\item 
A pixel that returns true for the following conditional statement is determined to be cloud obscured, therefore is not analysed further. 
\begin{align*}
    \rho_{0.65} + \rho_{0.86} > 1.2 \ &\text{or}\\
    T_{12} < 265  \ &\text{or}\\
     (\rho_{0.65} + \rho_{0.86} > 0.7 \ \text{and}\ T_{12} < 285 )\  &\text{or} \\
      (\rho_{0.86} > 0.25 \ \text{and}\ T_{12} < 300 )
\end{align*}
\item
If a pixel has a solar zenith angle that is less than 85 degrees, it is assumed to be a daytime pixel. A daytime pixel is considered a fire pixel if the following conditional statement is true.
$$ T_4 > 310 \ \text{and} \ (T_4 - T_{11}) > 10 \ \text{and} \ \rho_{0.86} < 0.35$$
\item 
If a pixel has a solar zenith angle that is greater than 85 degrees, it is assumed to be a nighttime pixel. A nighttime pixel is considered a fire pixel if the following conditional statement is true. 
$$ T_4 > 305 \ \text{and} \ (T_4 - T_{11}) > 10 $$
\item 
All other pixels are considered to be non-fire pixels. 
\end{enumerate}

\subsection{Implementing Edge Computing to Provide Fast and Efficient Fire Detection}
Computing speeds onboard satellites are much lower than ground-based computing speeds.  One of the fastest processors that are considered to be space-grade at the moment is the RAD750, which has a maximum CPU clock rate of only $200$ MHz \cite{bae}. Hence running the algorithm on a single satellite is inefficient. Thus, here we develop an edge-computing based method that uses several satellites for onboard bushfire detection. To do so, we require inter-satellite communication, which can be achieved by a geostationary satellite \cite{ec2}. With its high altitude, a geostationary satellite can maintain communications with a large number of LEO satellites in the constellation at one time.

\subsubsection{Choosing a Geostationary Orbit for Inter-satellite Communications}
Here, we choose our inter-satellite link to be a satellite in orbit similar to Optus D3. As this is a geostationary orbit it has a constant view of Australia at all points in time, allowing fire detections to be transmitted directly to Australia without delay \cite{D3}. The Keplerian coordinates of this satellite orbit are given in Table \ref{t3}.

\begin{table}[H] 
\caption{Keplerian coordinates of the geostationary satellite selected for inter-satellite communications}
\label{t3}
%%% \tablesize{} %% You can specify the fontsize here, e.g., \tablesize{\footnotesize}. If commented out \small will be used.
\begin{tabular}{cc}
\centering 
\textbf{Parameter}	& \textbf{Value}\\

Semi major axis ($a$)     & $ \SI{42165}{\kilo\meter}$ \\ Eccentricity ($e$)     & 0.0002541 \\
Inclination ($i$)  &  $\SI{0.0116}{\degree}$ \\ 
RAAN ($\Omega$) & $\SI{48.4858}{\degree}$ \\
Argument of perigee ($\omega$) & $\SI{135.8460}{\degree}$\\
Mean anomaly ($MA$) &$\SI{294.4219}{\degree}$ \\
\bottomrule
\end{tabular}
\end{table}

\subsubsection{Determining the LEO Satellite Visibility to the GEO Satellite}
At any given point in time, we need to know what LEO satellites of the constellation are visible to the GEO satellite to know what satellites can be allocated processing tasks. 
\par
\begin{figure}[h!]
    \centering
    \includegraphics[width=12 cm]{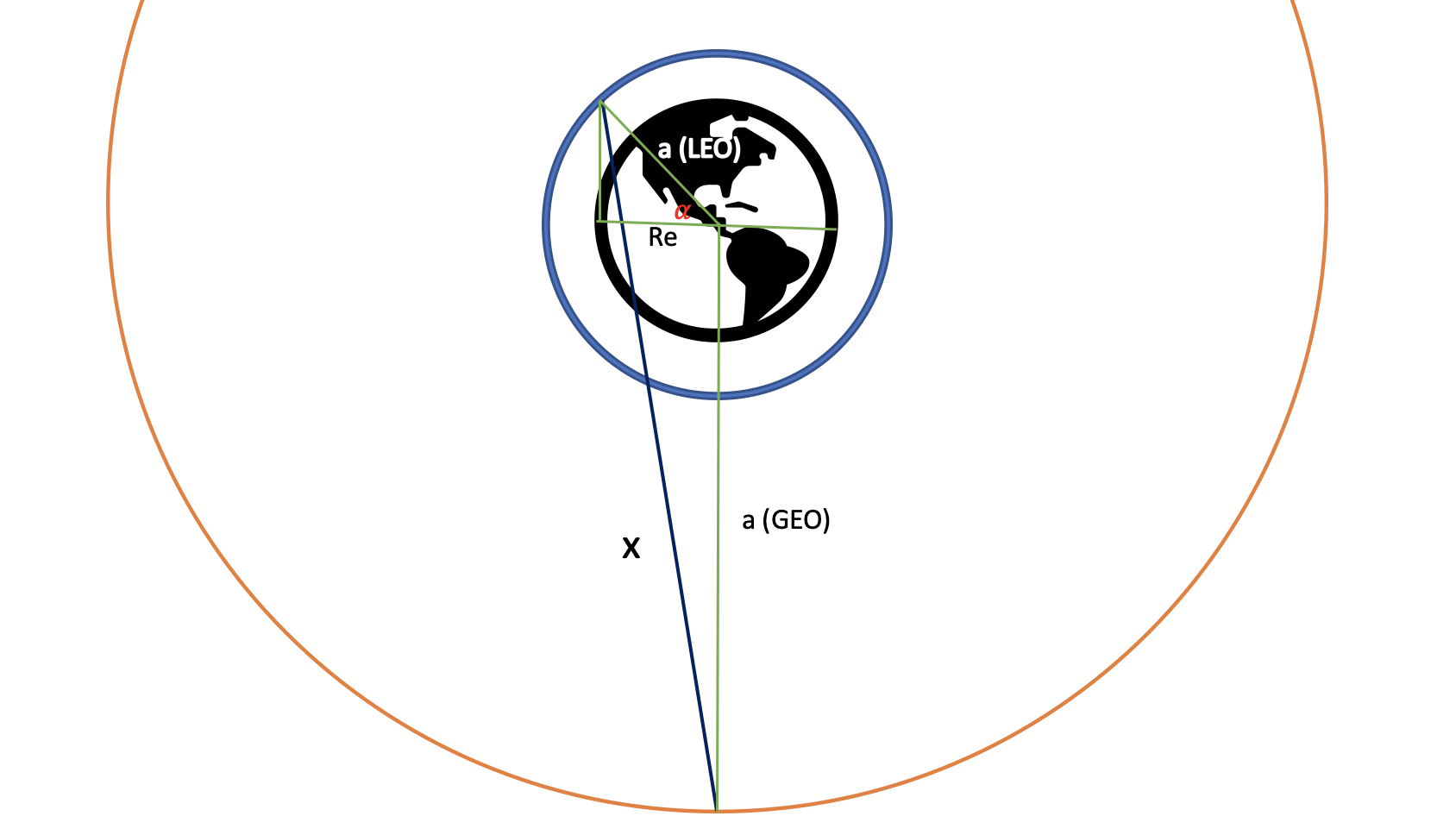}
    \caption{Visibility determination between LEO satellites and the GEO satellite. The LEO satellite orbit is shown in blue, while the GEO satellite orbit is shown in orange. x denotes the closest distance between the LEO and GEO satellites at a given time t.  }
    \label{vis}
\end{figure}

 
According to Figure \ref{vis}, when the distance between an LEO satellite of the constellation and the GEO satellite is greater than x, the satellites would not be visible to each other. From Figure \ref{vis}:
\begin{equation}
    \alpha = \cos^{-1}(\frac{R_e}{a(LEO)})
\end{equation}

Where $R_e$ is the radius of Earth and $ a(LEO)$ is the semi-major axis of the LEO orbit in question. We assume that both LEO and GEO orbits are circular to make our calculations simpler. This is a valid assumption as the eccentricities of all orbits involved are very low. Now by cosine rule: 
\begin{equation}
    x = \sqrt{ {\text{a(LEO)}}^2 + \text{a(GEO)}^2 - 2 \text{a(LEO)}\text{a(GEO)} \cos(90 + \alpha)}
    \label{eq1}
\end{equation}
If the distance between any LEO satellite and the GEO satellite is less than x, the satellites are visible to each other. Now we can use this to determine which LEO satellites are visible to the geostationary satellite so that inter-satellite communication links can be established between them.

\subsubsection{Steps of Our Edge Computing Algorithm}
\begin{figure}[h!]
    \centering
    \includegraphics[width=13 cm]{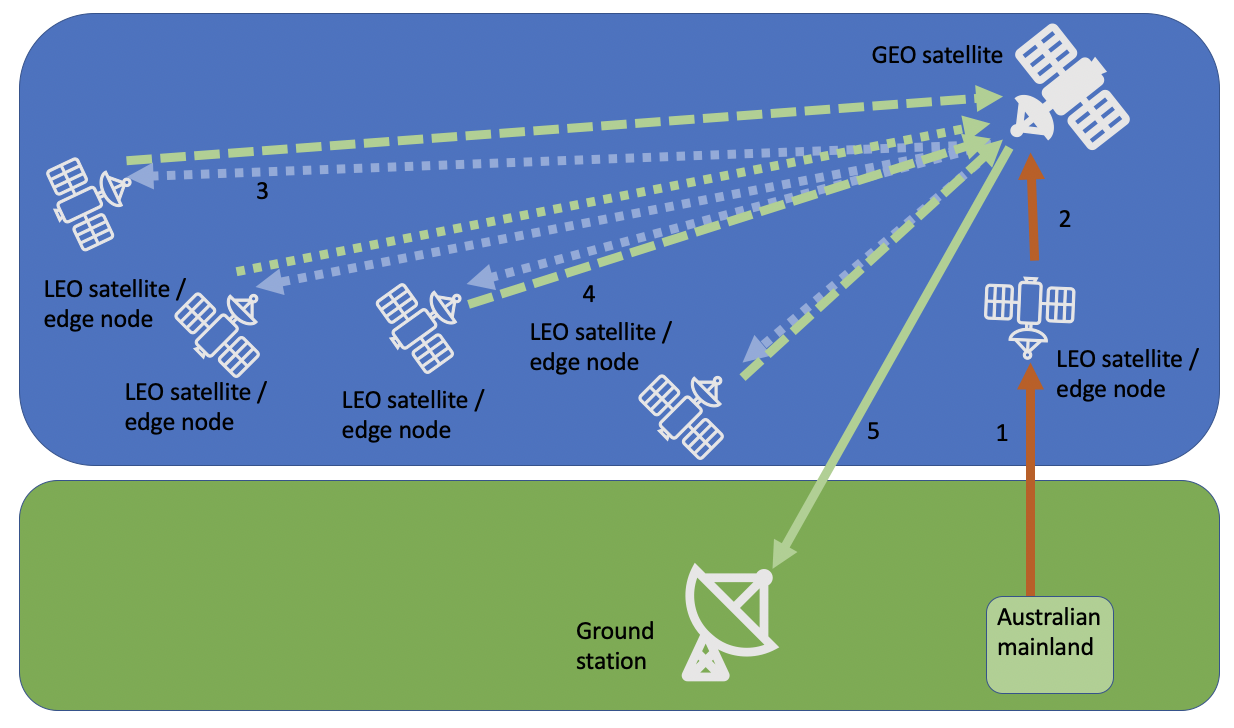}
    \caption{Inter-satellite comunication links established for the edge computing process }
    \label{pic}
\end{figure}

The following actions must take place in the given order to successfully distribute the image data among the constellation satellites, compile the results together, and transmit the outcome back to Earth. These steps are also illustrated in Figure \ref{pic}.

\begin{enumerate}
    \item A satellite above Australia captures an image extracts the reflectance,  radiance, latitude and longitude data for each pixel in the image and segments the image data to n parts.
    \item The satellite sends (n-1) parts to the GEO satellite. 
    \item The GEO satellite receives the (n-1) segments and identifies which LEO satellites are visible to it using the calculation described in the previous section, and sends each segment to a different visible LEO satellite in the network.
    \item Each satellite receives its segment of data and processes it using the Collection 6 algorithm described previously.  If a bushfire pixel is detected within the segment of data, the satellite relays the pixel's latitude and longitude back to the GEO satellite. 
    \item GEO satellite collects all bushfire coordinates it receives and sends the compilation down to an Australian groundstation. If no bushfire pixels are present, nothing is relayed to the ground, saving time and preventing redundant communications. 
\end{enumerate}
Each LEO satellite involved works as an edge node in this process, processing the data it receives through the GEO satellite. The GEO satellite works as the central edge node, allocating tasks to other edge nodes/satellites. What makes this an application of edge computing is the fact that none of the processing happens on the ground. Data is processed on the data receiver itself (and on other edge nodes).

%%%%%%%%%%%%%%%%%%%%%%%%%%%%%%%%%%%%%%%%%%
\section{Results}

In this section, we provide the results obtained using the methodology discussed in the previous section. The final constellation structure, the results and performance of the onboard image processing algorithm, as well as the reduction in detection time due to the use of the edge computing algorithm is discussed.

\subsection{Final Constellation Structure}
The final constellation structure derived through NSGA-II optimisation has 100\% coverage over the area of interest and 75.4\% coverage over the middle of Australia.  Its parameters are given in Table \ref{t6}.

\begin{table}[H] 
\caption{Final constellation orbits}
\label{t6}
%%% \tablesize{} %% You can specify the fontsize here, e.g., \tablesize{\footnotesize}. If commented out \small will be used.
\begin{tabular}{cc}
\toprule
\textbf{Parameter}	& \textbf{Value}\\
\midrule
Semi major axis ($a$)                   &$\SI{7334.9}{\kilo\meter}$\\
Altitude (h) & $\SI{963.9}{\kilo\meter}$\\ 
Eccentricity ($e$)                     & $ 0.04$ \\
Inclination ($i$)                     & $\SI{141.39}{\degree}$ \\
Number of orbital planes (P) & $ 95$ \\
Relative spacing between satellites in adjacent planes (F) & $9$ \\
Number of satellites per plane (n)                 &  $42$\\
Total number of satellites (N) & 3990\\ 
\bottomrule
\end{tabular}
\end{table}

\begin{figure}[h!]
   \centering
  \includegraphics[width=0.9\textwidth]{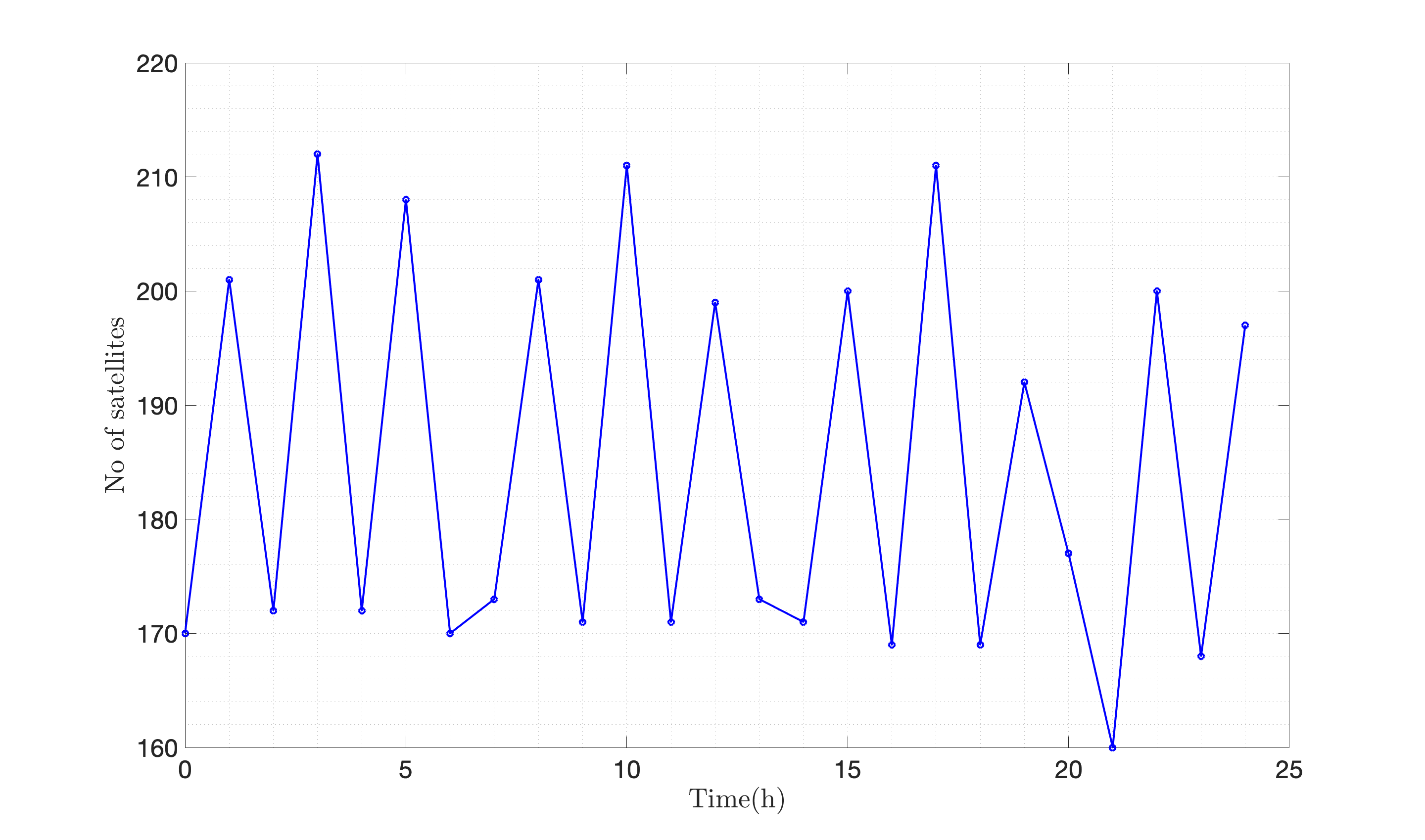}
  \caption{Number of satellite above Australia at each hour of the day}\label{3990vissat}
\end{figure}

\begin{figure}[h!]
     \centering
  \includegraphics[width=0.8\textwidth]{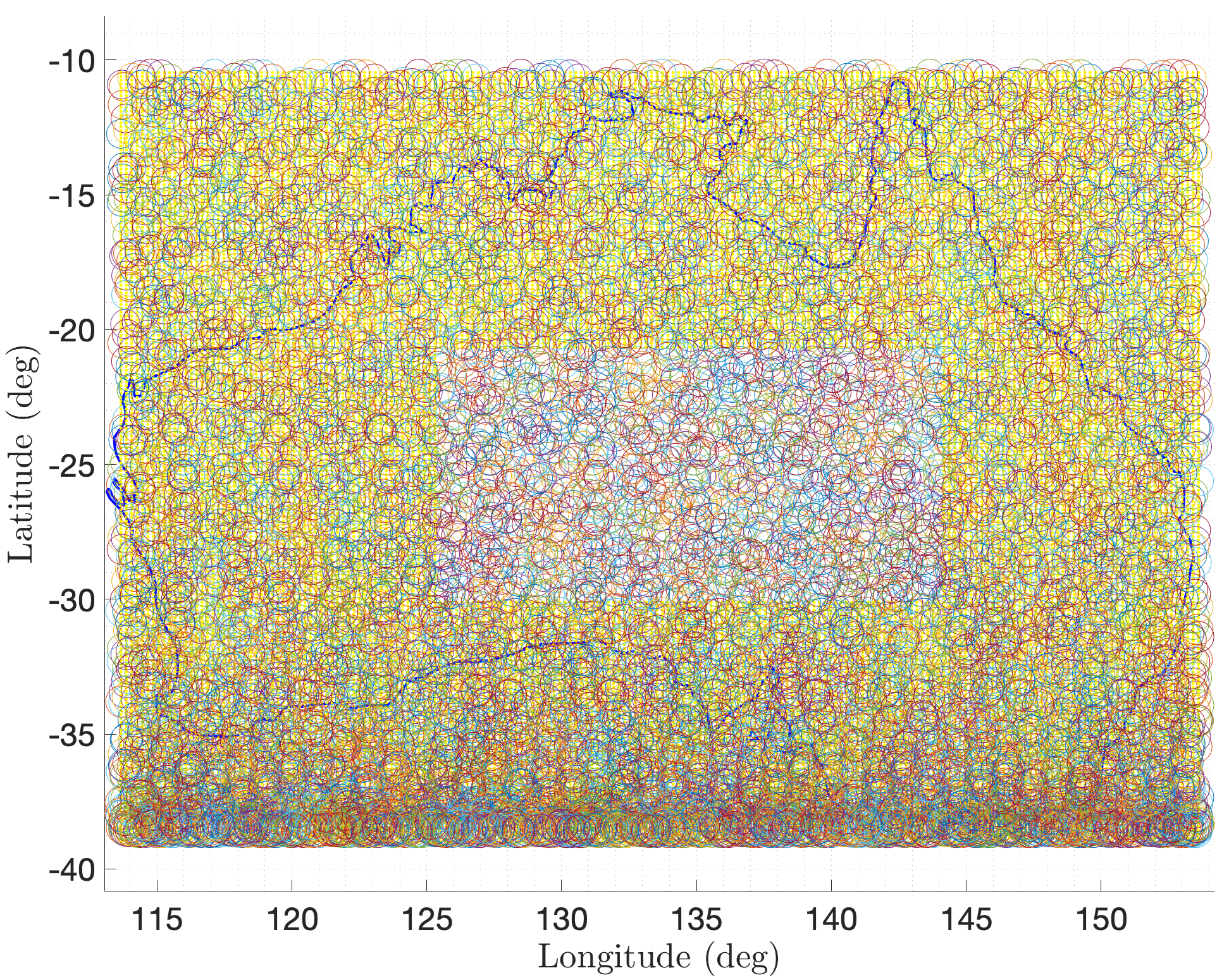}
  \caption{Coverage of the constellation (Area of interest is highlighted in yellow)}\label{3990aussie}
\end{figure}

\begin{figure}[h!]
    \centering
  \includegraphics[width=0.9\textwidth]{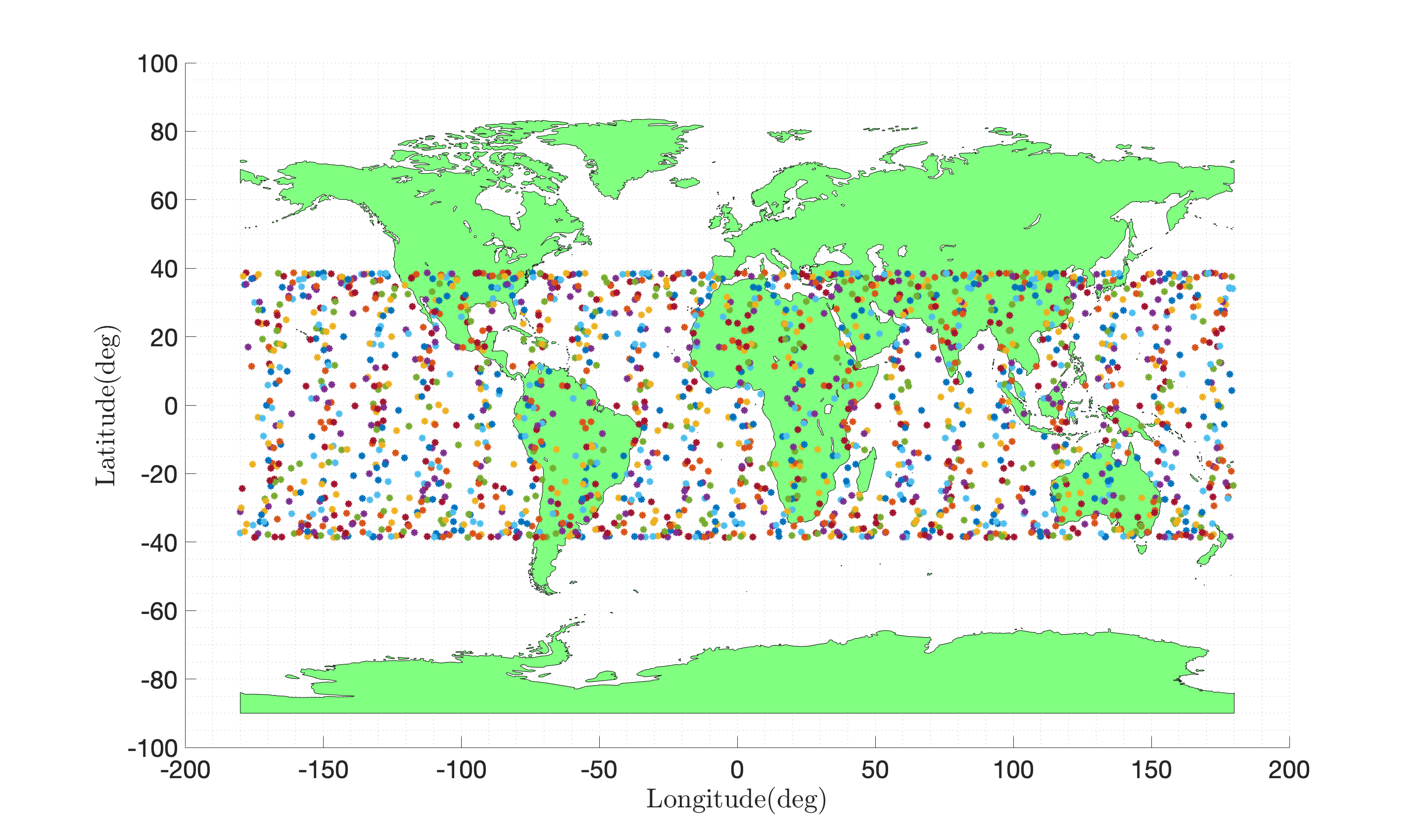}
  \caption{Ground track of the constellation}\label{3990gtrack}
\end{figure}

Figure \ref{3990vissat} shows how the number of satellites over Australia varies throughout the day for this constellation, while Figure \ref{3990aussie} shows its coverage over Australia when only the minimum number of satellites are visible to the area of interest. Figure \ref{3990gtrack} shows the ground track of 80 randomly chosen satellites of the constellation.  \par 
This constellation generates a new image of the area of interest once every 6 minutes and can generate a new image of central Australia once every 10 minutes.

\begin{figure}[H]
\centering 
\begin{tikzpicture}[thick]
\coordinate (O) at (0,0);
\coordinate (A) at (2,0);
\coordinate (B) at (0,4);
\coordinate (C) at (-2,0);
\draw (O)--(A)--(B)--(C)--cycle;
\draw (O)--(B);

\tkzLabelSegment[below=2pt](C,A){\textit{Swath Width}}
\tkzLabelSegment[left=2pt](O,B){\textit{h}}
\tkzLabelSegment[above right=2pt](A,B){\textit{}}
\tkzMarkAngle[fill= orange,size=0.7cm,%
opacity=.4](O,B,A)
\tkzLabelAngle[pos = 0.5](O,B,A){$\theta$}
\tkzMarkAngle[fill= orange,size=0.7cm,%
opacity=.4](C,B,O)
\tkzLabelAngle[pos = 0.5](C,B,O){$\theta$}
\end{tikzpicture}
\caption{Swath width calculation using the field of view angle, given by $2 \theta$}
\label{swath}
\end{figure}
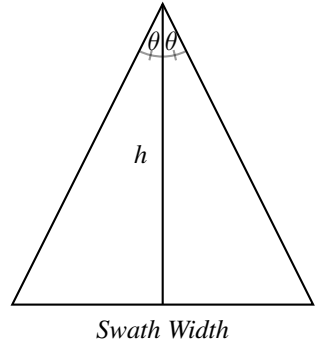

The Swath width of a satellite in the constellation can be calculated as follows using Figure \ref{swath}. 
$$ \text{Swath Width} = 2 h \tan{\theta} = 2\times  963.9\tan{6.99} = \SI{236.4}{\kilo\meter}$$
As this is greater than 44.4 km, the coverage criterion discussed in section \ref{cover} is fulfilled.

\subsection{Accuracy of the Onboard fire detection algorithm}\label{fire}
\begin{figure}{h!}
    \centering
    \includegraphics[width=0.9\textwidth]{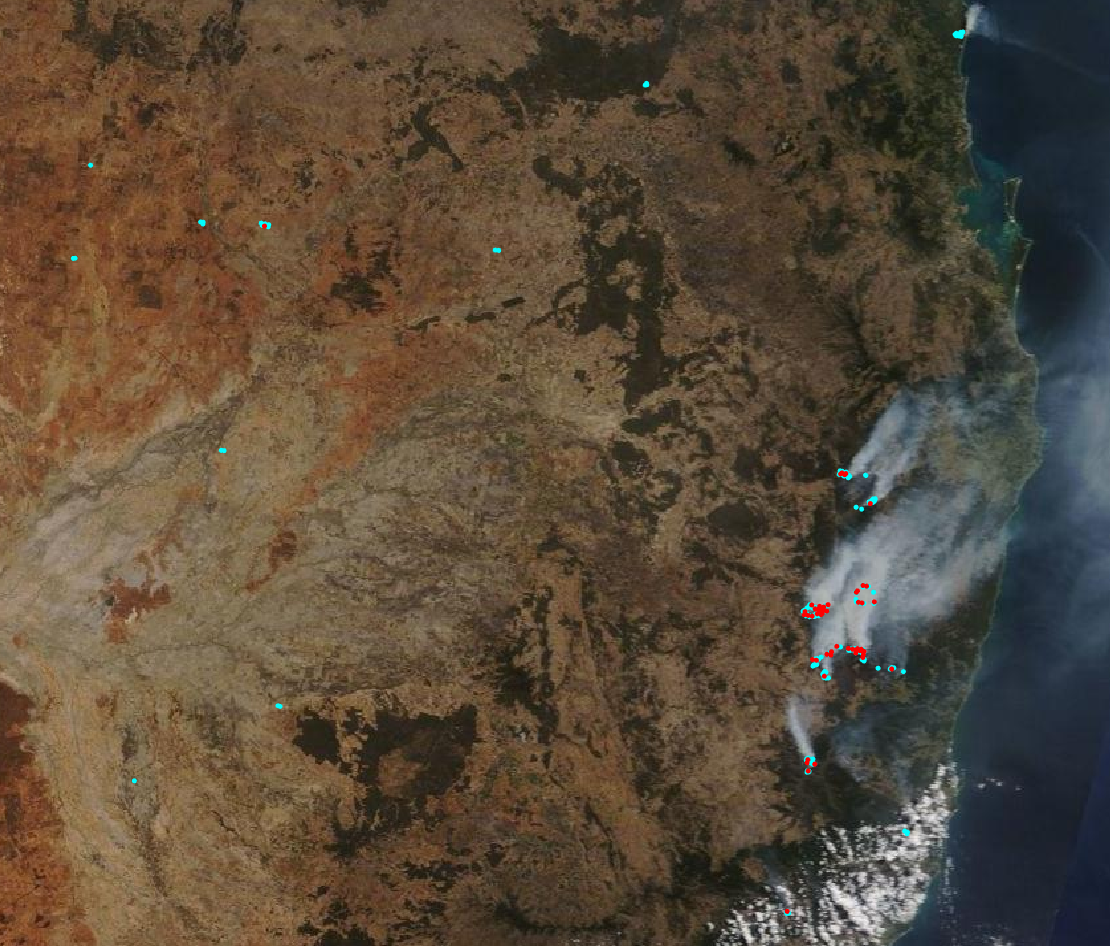}
    \caption{New South Wales as seen by Terra on 13th September 2019 \cite{Edata}, with the potential bushfire locations denoted in red.}
    \label{ss1}
\end{figure}

\begin{figure}{h!}
    \centering
    \includegraphics[width=0.9\textwidth]{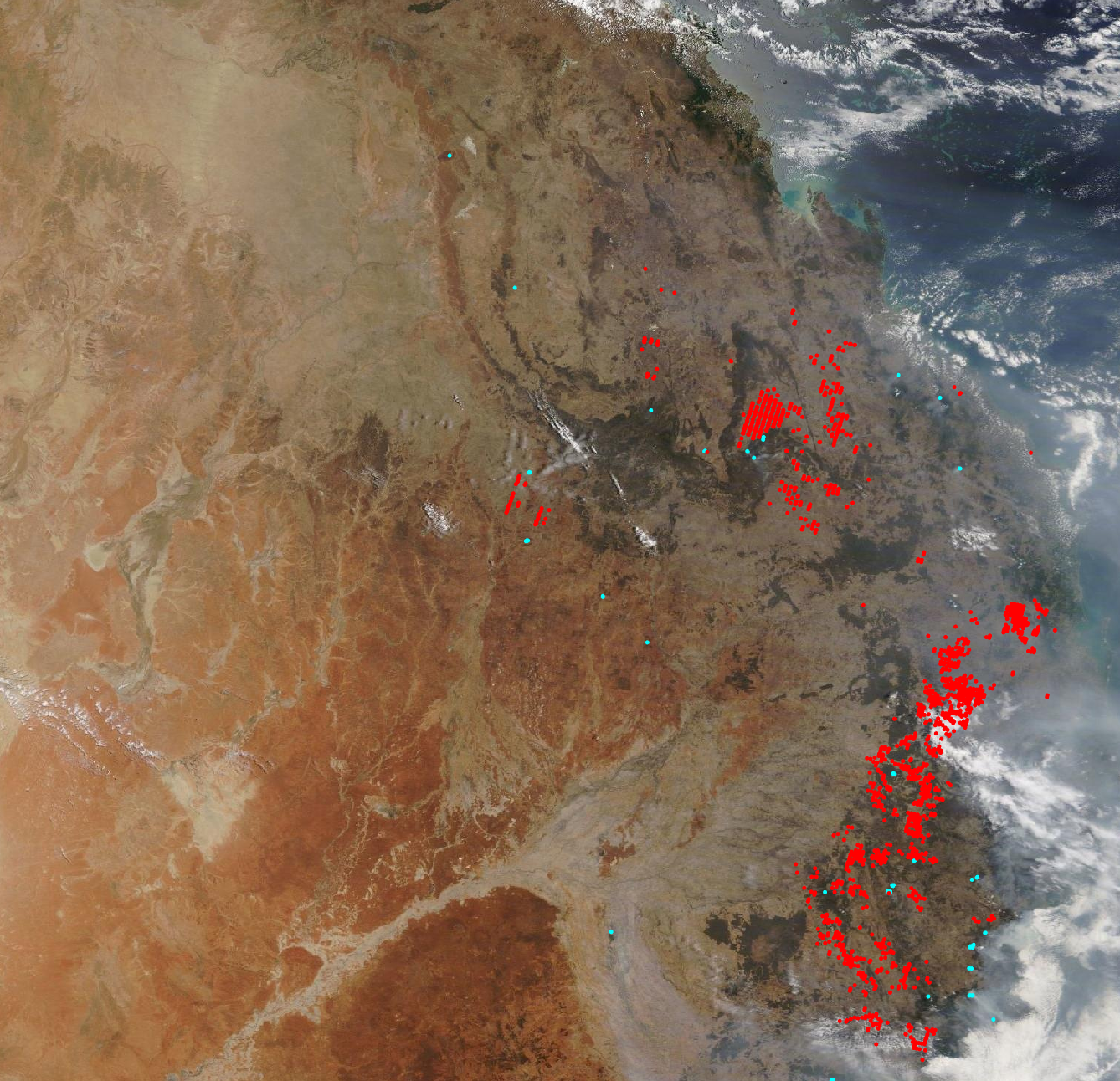}
    \caption{New South Wales as seen by Terra on 11th December 2019 \cite{Edata}, with the bushfires detected using the algorithm denoted in red.}
    \label{ss2}
\end{figure}

\begin{figure}{h!}
    \centering
    \includegraphics[width=0.9\textwidth]{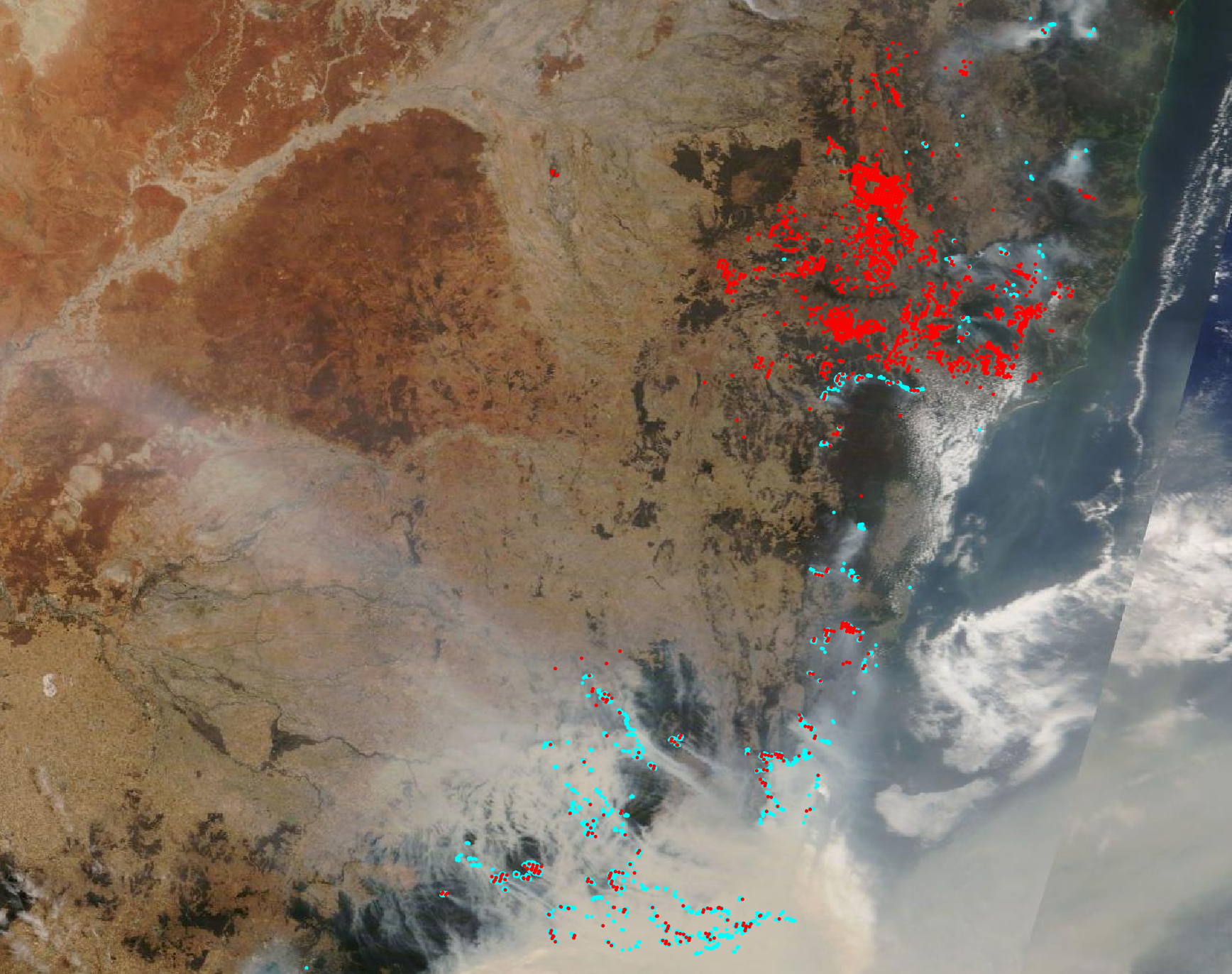}
    \caption{New South Wales as seen by Terra on 3rd January 2020 \cite{Edata}, with bushfire detections shown in red. Confirmed bushfires are shown in blue}
    \label{ss3}
\end{figure}

In this section, we provide the results from the image processing algorithm discussed in the methodology section. Figures \ref{ss1}, \ref{ss2} and \ref{ss3} show the fires detected. The red spots indicate fires detected using our algorithm, while the blue spots indicate known fires from NASA's Fire Information for Resource Management System (FIRMS) database \cite{firm}. \par

\begin{table}[H] 
\caption{Accuracy of the processed Images}
\label{t4}
%%% \tablesize{} %% You can specify the fontsize here, e.g., \tablesize{\footnotesize}. If commented out \small will be used.
\begin{tabular}{cc}
\toprule
\textbf{Image}	& \textbf{Accuracy}\\
\midrule
Image \ref{ss1}    &  $95.4\%$ \\
Image \ref{ss2}  &  $94.5\%$  \\ 
Image \ref{ss3}  &$59.3\%$\\
Average Accuracy  &$83.1\%$\\
\bottomrule
\end{tabular}
\end{table}

 The detection accuracy of the algorithm is summarised in Table \ref{t4}. The inaccuracies can be attributed to the fact that NASA's records of confirmed fires contain data from Terra, Aqua and other satellites, whereas our algorithm has only processed image data from the satellite Terra. Furthermore, NASA's records show fires happening throughout the day, while the image data processed by our algorithm is only relevant to a specific time of that day \cite{firm}.

\subsection{Results of the Edge Computing Algorithm}
In this section, we discuss the results obtained in the innovative process of using edge computing to allocate tasks of image processing among constellation satellites. 
Using Equation \ref{eq1}, we find that at a given point in time, at least 95\% of the LEO satellites (i.e. at least 3793 satellites) in the constellation are visible to the GEO satellite. \par 
Now we calculate the time taken to complete the five steps mentioned in the previous section, to see if our algorithm can provide improved fire detection time compared to current bushfire detection techniques.

\begin{figure}[h!]
    \centering
    \includegraphics[width=13.5 cm]{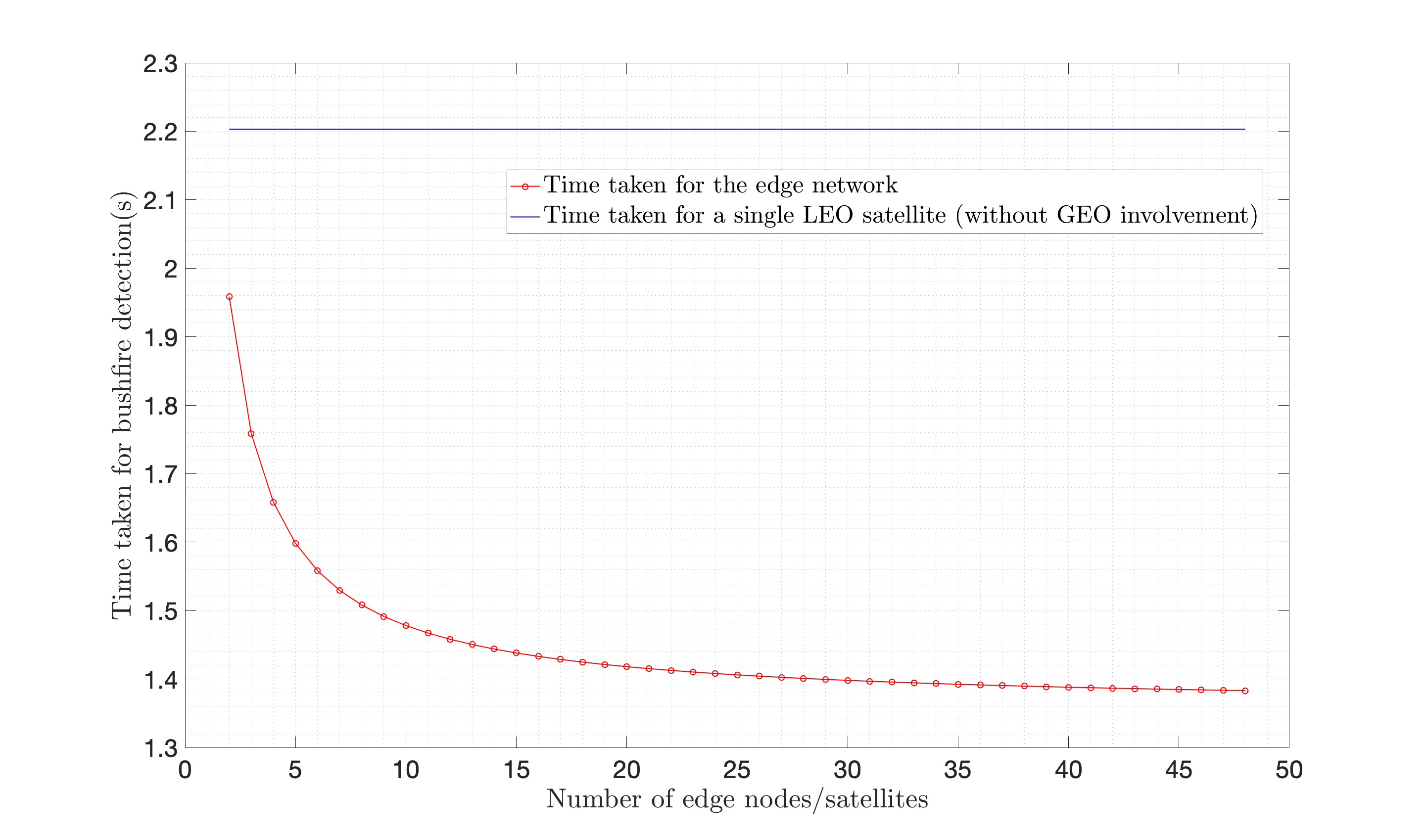}
    \caption{The relationship between the time taken to detect bushfires and the number of edge nodes used}
    \label{finale}
\end{figure}

\begin{itemize}
\item $t_1$ - Time is taken by the onboard sensor to collect sufficient photons from the ground to achieve the required resolution. This is known as the dwell time \cite{dt}. It can be calculated using Equation \ref{dwelltime}. 
\begin{equation}
    t_1 = \text{Dwell Time} = \frac{\text{Down track pixel size}/\text{Orbital velocity}}{\text{Cross-track line width }/ \text{ Cross-track pixel size}}
    \label{dwelltime}
\end{equation}
Our satellites are to have similar sensor capabilities to the RapidEye satellites as mentioned before. Thus, it is sensible to use a similar Multispectral push broom sensor for our satellites with a down track pixel size of 5 m. We shall calculate the lowest orbital velocity so that we can get an upper estimate of $t_1$. For this, we use the largest radial distance possible under the limits imposed on the chromosome used for the optimisation process.  $ r^{sat}_{min} = \SI{7375}{\kilo\meter}$
$$ \text{Max orbital velocity} = \sqrt{\frac{G M_{Earth}}{{r^{sat}_{min}}}}=   \sqrt{ \frac{\SI{6.673e-11}{} (\SI{5.98e24}{})}{\SI{7375e3}{}}} = \SI{7355}{\meter\per\second}  $$
For a push broom sensor, $\text{Cross-track line width }/ \text{ Cross-track pixel size}= 1$  \cite{dt}. We assume an additional 1 second is taken for segmenting the data into $n$ parts. Thus:  
$$ t_1 =\frac{5/7355}{1} + 1  = \SI{1.0006797}{\second}$$
\item $t_2$ - This is the time taken to transmit (n-1) segments of data from the LEO. It is sensible to assume that the data rate of our satellites would also be close to the RapidEye's data rate, which is 80 Mbit/s \cite{RE}.  For both RapidEye and MODIS, a single output file of data contains 559.3 kbit of data.
Thus: 
    $$ t_2 = \text{Time taken to transmit data from LEO Satellite} =  \frac{\text{Data amount}}{\text{LEO data rate}} = \SI{0.00699}{\second}$$
  \item $t_3$ - This is the time taken for the signal to travel from the LEO satellite to the GEO satellite.
    $$ t_3 = \text{Time taken for the data to travel from LEO to GEO} =  \frac{\text{Distance from LEO to GEO}}{ \text{Speed of light}} $$
The GEO satellite selected is at a semi major axis of $\SI{42165}{\kilo\meter}$. The LEO satellites are at a semi major axis of $\SI{7375}{\kilo\meter}$. Thus the distance between them is $\SI{34790}{\kilo\meter}$. Thus:$$t_3 = \SI{0.1160}{\second}$$
 
  \item$t_4$ -  This is the time taken for the  GEO satellite to receive the data. For a satellite like Optus D3, the communication bandwidth is around 54 MHz \cite{D3}. Assuming our satellite has the same bandwidth, we can determine the maximum possible data rate by Shannon's theorem \cite{shannon}.  
    $$ C = B \log_2(1 + SNR)$$
    Here, C is the maximum data rate, B is the bandwidth, and SNR is the signal to noise ratio. The typical SNR for RapidEye satellite data (which is now going to the GEO satellite) is around 210 \cite{modis1}. This gives us a data rate of :
     $$ C = B \log_2(1 + SNR) = 54\times 10^6\log_2(1 + 210) = \SI{416939356}{}\text{bps} \approx \SI{416.9}{}\text{Mbps} $$
    $$ t_4 = \text{Time taken for the GEO satellite to receive data} 
    $$$$ t_4 =  \frac{\text{Data amount}}{\text{GEO data rate}} = \frac{559300}{416939356} = \SI{0.00134}{\second}$$
    
   \item $t_5$ - This is the time taken for the GEO satellite to redistribute the data to other edge nodes for processing. We assume that the same amount of time is taken for the GEO satellite to redistribute the data to other edge nodes for processing. 
         $$  t_5 =  \text{Time taken for the GEO satellite to send away data} = t_4 = \SI{0.00134}{\second} $$
         \item $t_6$ -  This is the time taken for the signals to travel from the GEO satellite to the chosen LEO satellites.We can assume this to be same as $t_3$. 
       $$ t_6 = \text{Time taken for the data to travel from GEO to LEO} = t_3 = \SI{0.1160}{\second} $$

    \item$t_7$ -  Now, each LEO satellite takes time to receive the data from the GEO, and the time taken depends on the amount of data received. Assuming we have n edge nodes available, we get the following. 
           $$ t_7 =\text{Time taken for the LEO satellites to receive data} =   \frac{\text{Data amount}}{\text{LEO data rate}} = \frac{\frac{559300}{n}}{80\times 10^6} $$
\item$t_8$ - This is the time taken for processing the image data. 
To get an accurate estimate of how much time our algorithm would take to run onboard a satellite, we convert the bushfire detection code to C, compile it and then count the number of lines in the generated assembly (.asm) file. Assuming that a single clock cycle processes a single line of assembly code, we can roughly estimate the time that the algorithm would take to run onboard a satellite \cite{theiss}. \par 
Once the C code was compiled, it was noted that the generated assembly code had  \textbf{241878560} lines. The fastest processor that is considered to be space-grade at the moment is the RAD750, which has a maximum CPU clock rate of only 200 MHz \cite{bae}.
Thus, the time taken to run this algorithm can be calculated as follows: 
$$ t_8 = \text{Processing time} = \frac{1}{n} \frac{ \text{Number of lines in the assembly code}}{\text{Clock rate}} = \frac{1}{n}\frac{241878560}{200\times 10^6}$$
          \item $t_9$ - The  LEO satellites need time to send the coordinates of bushfire pixels back to the GEO satellite. We will assume here that the data sent away is the same size as the data received (It is likely less than this).
               $$ t_9 = \text{Time taken for the LEO satellites to send data} =  \frac{\text{Data amount}}{\text{LEO data rate}} =\frac{\frac{559300}{n}}{80\times 10^6} $$
               
\item $t_{10}$ -      Time is also taken for the signal to travel from LEO to GEO, we can again assume this to be the same as $t_3$. 
                 $$  t_{10} = \text{Time taken for the data to travel from LEO to GEO} = t_3 $$
    \item $t_{11}$ -    Time is taken for the GEO satellites to receive data and then -if there are bushfire pixels- send data to the ground.  Again assuming that the data received in this case is of the same size as the data received in $t_5$,
    $$ t_{11} = \text{Time taken for the GEO satellite to receive data} = t_5$$
    \item $t_{12}$ -    Time is taken for the GEO satellite to send the coordinates of fire pixels to the ground station. Again assuming size similarity, this time is equivalent to (and less than as the data size is likely smaller) $t_4$.
$$t_{12} = \text{Time taken for the GEO satellite to send data to ground}  = t_4$$
\item $t_{13}$ -   
         Lastly, time is also taken for the signal to travel to the ground. 
              $$ t_{13} = \text{Time taken for the data to travel from GEO to ground} $$$$ t_{13}=  \frac{\text{Distance from GEO to ground}}{\text{Speed of light}} = \frac{\SI{42165e3}{}}{c} = \SI{0.140}{\second} $$
         
\end{itemize}
By summing up the above times for different values of n, we can see how increasing the number of edge nodes (satellites) speeds up the process of bushfire detection. For the constellation structure described in this chapter, Figure \ref{finale} illustrates how the detection time reduces and plateaus out as the number of edge node increases. \par 
In Figure \ref{finale}, the blue line indicates the time taken to process an image onboard using just one satellite, without the use of the GEO satellite or other edge nodes. It is calculated by adding up the following times. 
\begin{itemize}
    \item $t_1'$ - This is once again the dwell time and the time taken to extract the reflectance and radiance data from an image. 
    Thus, $$t_1' = t_1 = \SI{1.0006796}{\second}$$
    \item  $t_2'$ - This is the time taken to process the data using just one satellite. 
    Thus, $$ t_2'  = \frac{\text{Number of lines in the assembly code}}{\text{Clock rate}} =  \SI{1.209}{\second}$$
    \item $t_3'$ - This is the time taken for the LEO satellite to send the data to the ground if a bushfire is detected. Both RapidEye and Modis has output files of 559300 bytes of data per image, and the LEO data rate is again assumed to be the same as that of RapidEye \cite{RE}. 
    $$ t_3' =  \frac{\text{Data amount}}{\text{LEO data rate}} =\frac{559300}{80\times 10^6} = \SI{0.00699}{\second} $$
    \item $t_4'$ - This is the time taken for the signal to travel from the LEO satellite to the ground station. 
       $$ t_4' =  \frac{\text{Distance from LEO to Ground}}{ \text{Speed of light}} = \frac{\SI{963.9e3}{}}{\SI{3e8}{}} = \SI{0.00321}{\second} $$
\end{itemize}
Adding up the above times, we get the time taken for a single satellite to process an image on-orbit to be 2.22 seconds. \par 
As seen in Figure \ref{finale}, the reduction in time taken plateaus to 1.39 seconds around 35 edge nodes.  Thus, there is no point in using more than 35 satellites in this process, as having more communication links than necessary could drastically reduce the power and resources of the satellites involved.
% add a table that compares this to that of Terra and Aqua

%%%%%%%%%%%%%%%%%%%%%%%%%%%%%%%%%%%%%%%%%%
\section{Discussion}
The final results of this project are summarised in Table \ref{t5}.

\begin{table}[h!] 
\caption{Accuracy of the processed Images}
\label{t5}
\begin{tabular}{lc}
\toprule
\textbf{Parameter}	& \textbf{Value}\\
\midrule
Image resolution   & $\SI{5}{\meter}$ \\
Swath width per satellite   & $\SI{236.4}{\kilo\meter}$ \\
Coverage over area of interest  &  $100\%$  \\ 
Coverage over central Australia  &  $75.4\%$  \\
Revisit time over area of interest  &  $\SI{6}{\minute}$ \\ 
Revisit time over central Australia  &  $\SI{10}{\minute}$ \\
Average fire detection accuracy  &  $83.1\%$  \\
Estimated fire detection time  &  $\SI{1.39}{\second}$ \\
\bottomrule
\end{tabular}
\end{table}

\subsection{ The Optimised Satellite Constellation}
The idea of using NSGA-II for orbit optimisation was inspired by a study on the construction of an  LEO-MEO constellation by M. Asvial \cite{ten}. However, all orbits considered in this study are circular, and only one parameter-the cost of construction- is optimised. The study conducted by Tania Savitri on-orbit optimisation is more similar to this project, as the parameter optimised is coverage, although global \cite{eleven}.  Both studies only consider constellations with satellites at the same altitude and inclination, while our optimisation occupies a much larger range of orbits, including varying altitudes, inclinations as well as eccentricities. Thus, in many ways, our project is an extension of previously conducted orbit optimisation studies.  \par 

The RapidEye satellite sensor was chosen to be used onboard the satellites, allowing us to obtain a small pixel size and a large swath width. Our pixel size is far smaller than that of Aqua, Terra, Sentinel, Himawari-8 or Landsat satellites; they all have pixel sizes that are upwards of 250 m \cite{Terra} \cite{him}. Hence, a bushfire must be at least 250 m wide for their detection. With the 5 m resolution, our constellation is capable of seeing bushfires that are 50 times smaller. \par 
Fires within the area of interest can be seen at any time by the constellation without significant delay. This constellation contains a large number of satellites, as monitoring the area of interest entirely at all times under 5 m resolution requires many satellites. However, considering that just the 2019-2020 bushfire season has estimated damages of 100 billion AUD \cite{bcost}, this would still be a feasible method of damage prevention. \par

\subsection{Onboard Image Processing Algorithm for Bushfire Detection}
 The image processing outcomes show that a simplified version of the Collection 6 algorithm can be successfully conducted on-board satellites to detect bushfires. This technique reduces the time and cost of bushfire detections, as we no longer have to downlink high-resolution images. \par 
As the FIRMS data source used to calculate the data accuracy of the detection algorithm contains fires detected throughout the day, comparing this data to the algorithm outcomes does not reveal much about its accuracy. Hence, while the majority of bushfires detected by our algorithm is correct, several confirmed fires are not identified by the detection algorithm (See Figures \ref{ss1}, \ref{ss2} and \ref{ss3}). These may be fires happening at other times of the day. There are also red dots far away from the blue ones, indicating that the algorithm can make infrequent false detections. These inaccuracies could potentially be eliminated by reinforcing the thresholding used in the algorithm with techniques such as machine learning methods (e.g. random forests). 

\subsection{Edge Computing Algorithm Image Processing}
Satellite constellations and the concepts of edge computing have been combined for efficiency in recent years. Several studies have been conducted on the use of edge computing for dynamic task allocation in satellite networks. For instance, the paper by Yongming He involves the use of edge computing on a satellite constellation to make observations upon user requests \cite{ec1}. Another example is the paper by Feng Wang that looks into the use of edge computing to handle networking requests from users \cite{ec2}.  The edge computing algorithms involved in both these studies are much different from ours as they have different end goals in mind. Furthermore, there are no current studies conducted where computations are taking place on the edge nodes themselves.
This concept is unique to our project. \par 
Figure \ref{finale} shows that increasing the number of satellites/edge nodes decreases the overall time taken for bushfire detections. However, this decrease asymptotes out at around 35 satellites, to a value of 1.39 seconds. Thus, having more than 35 edge nodes does not affect the efficiency of bushfire detection.  \par 
The reason for this distribution of processing was the fact that using a typical processor onboard a satellite, a single satellite would take more than 2.2 seconds to detect bushfires out of a single image obtained. As seen in Figure \ref{finale}, the time taken for the edge network-based processing is much less compared to this.  Thus, this combination of force solves the issue of time taken for onboard processing. \par 
With the use of a GEO satellite as the inter-satellite link, there is a constant ability to downlink data to the groundstation without delay. If the GEO satellite were not available, the LEO satellites would have to store the detection data until a groundstation becomes visible to them.\par 
The 1.39 s detection time is far less than the 1.5 h required for Aqua, Terra, Sentinel and Landsat satellites working in combination \cite{Terra}. It is even less compared to the 10 min detection time of the Himawari-8 satellite, which can only see bushfires that are at least 500 m long \cite{him}. However, having to establish multiple communication links may drain the power of the satellites far quickly compared to a single ground-LEO communication link. Hence, before implementation, calculations need to be conducted to see if this process is also advantageous from a power perspective.

%%%%%%%%%%%%%%%%%%%%%%%%%%%%%%%%%%%%%%%%%%
\section{Conclusions}
This project focuses on the construction of a network of satellites for bushfire detection in Australia. The constellation orbits are obtained through optimising the coverage, the total number of satellites and revisit times of the satellites. NSGA-II is chosen for the optimisation process, as it is more efficient compared to other optimisation algorithms. 
The constellation developed has 100\% coverage over coastal Australia and 75\% coverage over central Australia. It has 3990 satellites and can obtain an image of the coastal areas once every 6 minutes and of central areas once every 10 minutes. The constellation is at an altitude of $\SI{963.9}{\kilo\meter}$, an eccentricity 0.04 and an inclination $\SI{141.39}{\degree}$. It has 95 orbital planes, 42 satellites per orbital plane and relative phase separation of 9. 
All satellites are to be equipped with a crossbeam multispectral sensor similar to the RapidEye satellites, allowing them to provide real-time images with 5 m resolution in the visible and infrared wavelengths. The network is structured so that there are at least 160 satellites over coastal Australia at any given time. 
The satellites are also to be equipped with an on-orbit bushfire detection algorithm and processors with 200 MHz or higher clock rates. Edge computing is used to develop a novel method of distributed image processing, conducted through a GEO satellite in a similar orbit to Optus D3. Once a satellite above Australia takes an image, it is partitioned and sent to the geostationary satellite, which sends the segments off to 34 other constellation satellites to be processed. Each satellite processes its data segment, and if a bushfire pixel is detected, sends the fire coordinates to the GEO satellite. All bushfire coordinates are then compiled in the GEO satellite and downlinked to Earth. The entire process takes 1.39 seconds, leading to bushfire detections that are faster than any currently available method.

\section{Acknowledgements}
Please note that a preprint of this article been submitted as preprint in the below link: https://www.preprints.org/manuscript/202103.0071/v1 

\bibliographystyle{unsrt}  


\end{document}